\relax 
\documentstyle[12pt]{article}
\renewcommand{\baselinestretch}{1.2}

\newtheorem{defin}{Definition}[section]
\newtheorem{theorem}[defin]{Theorem}
\newtheorem{conj}[defin]{Conjecture}
\newtheorem{prop}[defin]{Proposition}
\newtheorem{lemma}[defin]{Lemma}
\newtheorem{cor}[defin]{Corollary}

\newcommand{\proof}{\noindent{\it Proof: }}
\newcommand{\qed}{{\mbox{$\ \ \rule{.35em}{1.7ex}$}}}

\newcommand{\rank}{\mbox{{\rm rank}}}

\newcommand{\oI}{\overline I}

\newcommand{\cR}{\widehat{R}}
\newcommand{\rR}{R_{red}}
\newcommand{\st}{(\star )}
\newcommand{\x}{\mbox{$x_{1},\ldots,x_{v}$}}
\newcommand{\xx}{\mbox{$x_{2},\ldots,x_{d}$}}
\newcommand{\ux}{\underline {x}}
\newcommand{\nnu}{\underline {\nu}}
\def\floor#1{\lfloor{#1}\rfloor}
\def\beneath#1#2{\mbox{\raisebox{.3ex}{$ \atop {\displaystyle#1 \atop #2}$}}}
\def\phiV{\phi\hskip-.04em\raisebox{-.3ex}{$\atop V$}}
\def\LAA{Linear Artin Approximation\ }

\catcode`\@=11
\input amssymb.sty
\catcode`\@=\active

\begin{document}

\vskip 1cm
\begin{center}
{\bf INTEGRAL CLOSURE OF IDEALS IN EXCELLENT LOCAL RINGS}
\vskip -2mm
{\bf (A corrected version)}
\vskip 2mm
Donatella Delfino and Irena Swanson
\vskip 2mm
\end{center}

This paper was published in
{\it J.\ Algebra}, {\bf 187} (1997), 422-445.
We are grateful to Ray Heitmann for pointing out that Theorem 2.7
in the published version is wrong.
Fortunately,
the main results of the paper are still true.
We give new proofs here.

\section{Introduction}

In \cite{BS}  Brian\c{c}on and Skoda proved, using analytic methods,
that if $I$ is an ideal
in the convergent power series ring ${\Bbb C} \{ x_{1},\ldots,x_{n} \}$
then ${\overline {I^n}}$, the integral closure of $I^n$,
is contained in $I$.
Extensive work has been done in
the direction of proving ``Brian\c{c}on-Skoda type theorems", that is,
statements about
${\overline {I^{t}}}$ being contained in $(I^{t-k})^{\# }$,
where $k$ is a constant independent of $t$,
and $\#$ is a closure operation on ideals
(cf. \cite{LSBriancon}, \cite{LTBriancon}, \cite{ReesSally}, \cite{SwBS}).

In this paper we study the following related problem:
given an ideal $I$ of a Noetherian local ring $(R,m)$, find a 
``linear'' integer-valued function $f(n)$
such that
$\overline{I+m^n} \subseteq \oI + m^{f(n)}$
for all $n$
or for all sufficiently large $n$.

An element $x$ of $R$ is said to be in the {\it {integral closure}}
${\overline {J}} $ of an ideal $J$
if it satisfies a relation of the form
$x^{n}+\alpha _{1}x^{n-1}+\alpha _{2}x^{n-2}+\ldots +\alpha _{n}=0$,
with $\alpha _{t}\in J^{t}$ for all $t$.
%
%
%
%
%
We first observe that
if $(R,m)$ is a noetherian local ring
which is complete in the $m$-adic topology,
then there exists an integer-valued function $f(n)$,
with $\lim_{n\to \infty} f(n)=\infty$,
such that
${\overline {I+m^{n}}} \subseteq \oI +m^{f(n)}$.
To see this, we use the fact that
$\oI = \beneath{\cap}{V} \phiV \hskip-.5em^{-1}(IV)$,
where the intersection is over all discrete valuation domains $V$
which are $R$ algebras via $\phiV$
and whose maximal ideal contracts to $m$.
With this it is easily shown that
$\beneath{\cap}{n} {\overline {I+m^{n}}}=\oI $.
By Chevalley's theorem (\cite[p.\ 270]{ZSII})
(for $R/\oI $
and the descending sequence of ideals $\{\overline{I+m^n}\,/\,\oI\}_n$)
then
${\overline {I+m^{n}}}\subseteq \oI +m^{f(n)}$
for some function $f(n)$ such that $\lim_{n\to \infty} f(n)=\infty$.
Chevalley's theorem does not help in determining the order of growth of $f(n)$,
and, in fact,
it cannot because it takes into account only the topology determined
by a given descending sequence of ideals.
By also considering the algebraic properties of the sequence
we prove a stronger statement:

\bigskip

\noindent{\bf The main theorem.}
Let $(R,m)$ be an excellent local ring.
Let $I$ be an ideal of $R$.
Then there exists a positive integer $c$ such that
$$
\overline{I+m^n} \subseteq
\oI +m^{{\floor{{n/c}}}}
{\mbox{\rm \ \ for all \ }} n.
$$
\bigskip

If $I$ is $m$-primary, then the theorem holds trivially.
Indeed,
large powers of $m$ are contained in $I$,
so
$\overline{I+m^n} = \oI \subseteq \oI + m^n$ if $n >> 0$.

Rees proved in~\cite{Reesanalunram} that
if $J$ is an ideal of an analytically unramified local noetherian ring
then ${\overline {J^{n}}}\subseteq J^{n-k}$
for a constant $k$ and all $n\geq k$.
If $J=m$ then Rees' theorem yields a special case of the main theorem
of this paper.

In general we only know the existence of a positive integer $c$
for which the theorem is satisfied,
but we can give an explicit bound for some classes of ideals
(see Examples 4.1, 4.2, Proposition \ref{prop:monomial}).

Our proof of the main result would be greatly simplified
(see comments after Proposition \ref{prop:dvr})
if we could use the following conjecture,
known as the \LAA Theorem:

\begin{conj}
{\em {(\LAA Theorem)}}
Let $(R,m,k)$ be a complete local ring.
Suppose we have a system of (finitely many) equations
in $t$ variables over $R$
and
we know that the system has a solution ${\underline {Z}}$ modulo $m^{l}$.
Assume that
if $J$ is the ideal generated by the polynomials defining the equations,
then $J\cap R=(0)$.
Then the system has a true solution ${\underline {U}}$ such that
${\underline {Z}}-{\underline {U}}\in m^{\floor {l/c}}R^t$,
where $c$ is a constant independent of $l$.
\end{conj}

The \LAA theorem was announced by Spivakovsky in \cite{Slaa}.
However, there is no proof of the theorem,
in that generality, in the literature.
Lejeune-Jalabert and Hickel proved the case when the ring
is an isolated hypersurface singularity (cf. \cite{Ljlaa}, \cite{Hickel})
and gave an explicit bound for $c$.
We use Lejeune-Jalabert's bound in Example 4.12.
We get around using the general \LAA Theorem by
actually proving a special case of it in the process
(see comment after Proposition \ref{prop:dvr}).

\bigskip
We close the introduction by summarizing the structure of this paper.
Section 2
proves that it suffices to show that the main theorem
holds for principal ideals in complete integrally closed domains.
Section 3 is devoted to the proof of the main theorem for
principal ideals in rings as above
and also contains the proof of a special case of the \LAA Theorem.
Section 4 provides explicit bounds for the
constant $c$ in several cases of interest.

\bigskip
\bigskip
The authors are grateful to Melvin Hochster and Craig Huneke for
many conversations regarding this material.
We learned much from Professor Hochster's insight.

\section{Some reductions}

In this section we prove that it is sufficient to prove the main theorem
in case $(R,m)$ is a complete normal (that is, integrally closed)
local domain and $I$ is a principal ideal.
We do this in several steps.

We start with a lemma which may justify the modified notation
$\st$ below:

\begin{lemma}
\label{lemma:radical}
Let $R$ be a Noetherian ring and $I$ and $J$ arbitrary ideals.
Let $K$ be the radical of $J$.
Then there exists an integer valued $f(n)$ tending to infinity such that
$\overline{I+K^n} \subseteq \oI +K^{f(n)}$
for all $n$
if and only if
there exists an integer valued $g(n)$ tending to infinity such that
$\overline{I+J^n} \subseteq \oI +J^{g(n)}$
for all $n$.

Moreover,
$f$ grows linearly if and only if $g$ does.
\end{lemma}

\proof
As $R$ is Noetherian,
there exists an integer $k$ such that $K^k \subseteq J$.

First we assume the existence of $f$.
Then
$$
\overline{I+J^n} \subseteq
\overline{I+K^n} \subseteq
\oI +K^{f(n)} \subseteq
\oI +J^{\floor{f(n)/k}}.
$$

Now we assume the existence of $g$.
Then
$$
\overline{I+K^n} \subseteq
\overline{I+J^{\floor{n/k}}} \subseteq
\oI +J^{g(\floor{n/k})} \subseteq
\oI +K^{g(\floor{n/k})}.
$$

In both cases it is clear that linear growth of one function implies
the linear growth of the other.	
\qed

We now set up some notation to express this and other more general cases:
In a Noetherian ring $R$ with ideals $I$ and $J$
we consider the existence of a constant $c$ such that

$$
\overline{I+J^n} \subseteq \oI +J^{\floor{n/c}}
\eqno{\st }
$$
for all $n$.					
If there exists such a $c$,
we say that $\st $ holds in $R$ for $I$ and $J$.

With this,
the lemma implies that $\st $ holds in $(R,m)$ for $I$ and an $m$-primary
ideal
if and only if the main theorem holds in $R$.

We use notation $\st $ also in the following reduction to the principal ideal case:

\begin{prop}
\label{prop:principal}
If $\st $ holds in every excellent local ring $(R,m)$
for every principal ideal $I$ and for $m$,
then $\st $ holds in every excellent local ring $(R,m)$
for every ideal $I$ and for $m$.
\end{prop}

\proof
Let $I$ be an arbitrary ideal in an excellent local ring $(R,m)$.
We want to prove that $\st $ holds in $R$ for $I$ and $m$.

We first let $S$ be the extended Rees ring $R[It,t^{-1}]$,
where $t$ is an indeterminate over $R$.
Let $M$ be the maximal homogeneous ideal $mS+ItS+t^{-1}S$ of $S$.
As $S$ is a finitely generated $R$ algebra,
$S$ and $S_M$ (localization at $M$) are both excellent rings.
Moreover $IS \subseteq t^{-1}S$.
By assumption $\st $ holds in $S_M$ for $t^{-1}S_M$ and $MS_M$.
This means that there exists a positive integer $c$ such that
$$
\overline{t^{-1}S_M+M^nS_M} \subseteq \overline{t^{-1}S_M} + M^{\floor{n/c}}S_M
$$
for all $n$.

Now
let $x$ be in $\overline{I+m^{n}}$.
Then as $R \subseteq S_M$,
$x \in \overline{IS_M+M^nS_M} \subseteq \overline{t^{-1}S_M+M^nS_M} \subseteq
\overline{t^{-1}S_M} + M^{\floor{n/c}}S_M$.
To finish the proof it thus suffices to show
$\left(\overline{t^{-1}S_M} + M^{\floor{n/c}}S_M\right) \cap R =
\oI + m^{\floor{n/c}}$.
First observe that by
the one-to-one correspondence of $M$-primary ideals in $S$
and $MS_M$ primary ideals in $S_M$ we have
$\left(\overline{t^{-1}S_M} + M^{\floor{n/c}}S_M\right) \cap S =
\overline{t^{-1}S} + M^{\floor{n/c}}$.
Thus it suffices to prove that
$\left(\overline{t^{-1}S} + M^{\floor{n/c}}\right) \cap R = \oI + m^{\floor{n/c}}$.
But
$\overline{IS} + M^{\floor{n/c}} =
\overline{IS} + (mS+It)^{\floor{n/c}}$
is a graded ideal whose graded piece of degree~0 is
$\overline{I} + (m)^{\floor{n/c}}$,
which finishes the proof.
\qed
\medskip

Thus from now on we may assume that $I$ is a principal ideal.
The next goal
is to replace $R$ by a complete Noetherian local domain.
The first step is to pass to the completion of $R$.
This will ensure that all the relevant finitely generated $R$-algebras
which are reduced
have module-finite integral closures in their total rings of fractions.
We use this property on $\rR$,
the quotient of $R$ by the ideal $\sqrt{0}$ of the nilpotent elements.
By using the Artin-Rees lemma
we conclude that then the main theorem holds in $R$
if it holds in the integral closure $S$ of $\rR$.
This $S$ is a direct sum of domains,
each of the domains being of the form $S/P$ for some minimal prime $P$ of $S$.
We prove that if $\st $ holds in each $S/P$ for the image of $I$
and the maximal ideal of $S/P$,
then $\st $ also holds in $S$ for $IS$ and $mS$.

Since the image of a principal ideal in any algebra is still principal,
by the reductions above
we end up with a principal ideal in a complete local normal domain.
The rest of this section is just proving that we may make these reductions.

({\it Comment:}
If we try to go modulo all the minimal primes before normalizing,
we cannot conclude $\st $ for $R$ from knowing $\st $ in all domain quotients,
as $R$ need not be a direct sum of such domains.)

\begin{lemma}
\label{lemma:complete}
Let $(R,m)$ be an excellent local ring
and let $\cR$ be the $m$-adic completion of $R$.
If the theorem holds in $\cR$,
it also holds in $R$.
\end{lemma}

\proof
By assumption there exists a positive integer $c$ such that
$\overline{I\cR+m^n\cR} \subseteq \overline{I\cR} +m^{\floor{n/c}}\cR$
for all $n$.
Now let $x\in \overline{I+m^{n}}$.
Then
\begin{eqnarray*}
x \in \overline{I+m^{n}} \cR \cap R &\subseteq&
	\overline{(I+m^{n})\cR} \cap R \\
&\subseteq& \left( \overline{I\cR} +m^{\floor{n/c}}\cR\right) \cap R \\
&=& \left( \oI\cR +m^{\floor{n/c}}\cR\right) \cap R \\
&&
\mbox{\rm \hskip5em
(by \cite[Examples v, iv p.\ 800 and Lemma 2.4]{Lipmanlipsat})} \\
&=& \oI +m^{\floor{n/c}}.
\qed
\end{eqnarray*}

Lemma \ref{lemma:complete} is the only place where the excellence
of the ring is used.

\begin{lemma}
\label{lemma:reduced}
Let $(R,m)$ be a Noetherian local ring.
If the theorem holds in $\rR = R/\sqrt 0$,
then it also holds in $R$.
\end{lemma}

\proof
By assumption there exists a positive integer $c$ such that
$\overline{I\rR+m^n\rR} \subseteq \overline{I\rR} +m^{\floor{n/c}}\rR$
for all $n$.
Now we use the fact that for any ideal $I$ in $R$,
the integral closure of $I$ in $R$
is the same as
the preimage in $R$ of the integral of $I\rR$ in $\rR$.
This implies
$$
\overline{I+m^{n}} \subseteq
\mbox{preimage of } (\overline{I\rR} +m^{\floor{n/c}}\rR)
= \oI +m^{\floor{n/c}}.
\qed
$$

These two lemmas say that it is enough to prove the main theorem
for complete local reduced rings.
The next lemma will enable us to normalize such a ring:

\begin{lemma}
\label{lemma:modulefinite}
Let $(R,m)$ be a Noetherian local ring and $I$ an ideal in $R$.
If $(\star )$ holds in a module-finite extension $S$ of $R$
for $IS$ and $mS$,
then the main theorem holds in $R$.
\end{lemma}

\noindent
{\it Note}
that $S$ may not be local.

\proof
As $S$ is a finite $R$-module,
$\overline{IS} \cap R = \oI$.
Also, the inclusion $R/\oI \subseteq S/\overline{IS}$ is module-finite
so by the Artin-Rees Lemma there exists an integer $k$ such that
$(m^n S + \overline{IS})/\overline{IS} \cap R/\oI
\subseteq (m^{n-k} R+ \oI) /\oI$ for all $n$.

Now we use the assumption that there exists a positive integer $c$ such that
$$
\overline{IS +(mS)^{n}} \subseteq \overline {IS}+(mS)^{\floor{n/c}}
$$
for all $n$.
Thus
$$
\overline{I +m^n} \subseteq
\overline{IS +(mS)^n} \cap R \subseteq
(\overline{IS}+(mS)^{\floor{n/c}}) \cap R
$$
for all $n$.
We rewrite this modulo $\oI$ and $\overline{IS}$:
$$
\overline{I +m^n} R/\oI \subseteq
(m^{\floor{n/c}} S + \overline{IS})/\overline{IS} \cap R/\oI
\subseteq (m^{\floor{n/c}-k} + \oI )/\oI.
$$
Hence
$\overline{I +m^n} \subseteq \oI + m^{\floor{n/c}-k}$.
It is easy to show that there exists a positive integer $c'$
($c(k+1)$ will do)
such that
$\floor{n/c}-k \ge \floor{n/c'}$ for all $n$.
Thus the final version says that
$\overline{I +m^n} \subseteq \oI + m^{\floor{n/c'}}$
for all $n$.
\qed
\bigskip

We have now reduced to the following situation:
$S$ is the integral closure of a reduced complete local ring $(R,m)$
in its total field of fractions.
Thus $S$ is an integrally closed reduced Noetherian ring,
module-finite over $R$ and
with finitely many maximal prime ideals $P_1$, \ldots, $P_l$
all containing $mS$.
Also,
$S$ is complete in the $mS$-adic topology
and $S/mS$ has dimension zero.
It follows that $S$ is a direct sum of finitely many domains,
each domain being of the form $S/P$ for some minimal prime $P$ of $S$,
or better yet,
each domain being of the form $S_{P_i}$ for some $i$:
$$
S = \lim_{\longleftarrow} S/{m^k S} =
\lim_{\longleftarrow} S/{P_1^k} \times \cdots \times
\lim_{\longleftarrow} S/{P_l^k}
=
S_{P_1} \times \cdots \times S_{P_l}.
$$
Note that $mS_{P_i}$ is $P_i$-primary.
The following lemma will reduce the proof of the main theorem to these complete
normal local domains:

\begin{lemma}
\label{lemma:product}
Let $R = R_1 \times \cdots \times R_l$
be a direct sum of rings.
Let $I$ and $m$ be ideals in $R$.
If $(\star )$ holds in each $R_i$ for $IR_i$ and $mR_i$,
then it also holds in $R$ for $I$ and $m$.
\end{lemma}

\proof
We will use the fact that for any ideal $J$ in $R$,
$$
JR = JR_1 \times \cdots \times JR_l
\mbox{\ \ \ and \ \ \ }
\overline{JR} = \overline{JR_1} \times \cdots \times \overline{JR_l}.
$$

By assumption there exist positive integers $c_i$ such that
$$
\overline{IR_i+m^nR_i} \subseteq \overline{IR_i}+m^{\floor{n/c_i}}R_i
$$
for all $n$ and all $i = 1, \ldots, l$.
Let $c = \max\{c_1, \ldots, c_l\}$.
Thus the inclusions above also hold when each $c_i$ is replaced by $c$.
With this,
\begin{eqnarray*}
\overline{I+m^n} &=&
\overline{IR_1+m^nR_1} \times \cdots \times \overline{IR_l+m^nR_l} \\
&\subseteq&
(\overline{IR_1}+m^{\floor{n/c}}R_1)
\times \cdots \times
(\overline{IR_l}+m^{\floor{n/c}}R_l) \\
&=&
\overline{IR_1}
\times \cdots \times
\overline{IR_l}
\ +\ %
m^{\floor{n/c}}R_1
\times \cdots \times
m^{\floor{n/c}}R_l \\
&=&
\oI+m^{\floor{n/c}}.
\qed
\end{eqnarray*}

So we reached the main goal of Section~2:
we started with an arbitrary excellent local ring $(R,m)$
with an arbitrary ideal $I$.
By Proposition \ref{prop:principal} we may assume that $I$ is principal.
By Lemmas \ref{lemma:complete}
and \ref{lemma:reduced} we may assume that $R$ is complete in the $m$-adic
topology
and that it has no nonzero nilpotents.
By Lemmas~\ref{lemma:modulefinite}, \ref{lemma:product}
and~\ref{lemma:radical},
we may then assume that $R$ is a complete normal local domain.

With these reductions,
the proof of the main theorem for principal ideals $I$
in complete local normal domains
is given in the next section.

\section{Proof of the main theorem}

Before proving the main theorem for complete normal local domains,
we consider some special cases:

\begin{prop} [Rees]
\label{prop:selfrad}
Let $I$ be a radical ideal in a complete Noetherian local ring $(R,m)$.
Then there exists an integer $k$ such that
for all $n\geq k$,
$\overline{I+m^n} \subseteq I+m^{n-k}$.
\end{prop}

\proof
Let $(-)'$ denote going modulo $I = \sqrt{I} = \oI$.
Then $R'$ is a complete reduced local ring,
hence analytically unramified.
Thus by~\cite[Theorem~1.4]{Reesanalunram}
there exists an integer $k$ such that
$\overline {(m')^n} \subseteq (m')^{n-k}$.
It follows that:
$$
(\overline {I+m^{n}})' \subseteq
\overline {I'+(m')^n} = \overline {(m')^n}
\subseteq (m')^{n-k},
$$
hence
$\overline {I+m^{n}} \subseteq m^{n-k} + I$.  \qed

One may naively think that the inclusion
$\overline {I+m^{n}_{R}} \subseteq \oI+m^{n-k}_{R}$
for some constant $k$ and all $n$ sufficiently large holds in general.
However, this is false.
Indeed, let $K$ be a field, $x, y$ variables over $K$, and $R=K[[x,y]]$.
Then $(xy^{n/2})^{2}=x^{2}y^{n}\in (x^{2}R+m^{n}_{R})^{2}$
for all even positive integers $n$,
so $xy^{n/2}\in \overline {x^{2}R+m^{n}_{R}}$.
However, there is no constant $k$ such that
$xy^{n/2}\in (x^{2},xy^{n-k-1},y^{n-k})R$
for all $n>>0$.

\begin{prop}
\label{prop:monomial}
Let $(D,\mu)$ be an analytically unramified local ring
and let $x$ be a variable over $D$.
Let $R=D[[x]]$ and $m = \mu R + (x)R$ the maximal ideal of $R$.
Then $\overline {x^tR+ m^n} \subseteq x^tR+m^{\floor{n/c}}$
for some positive integer $c$ and for all $n$.

Moreover,
if the powers of $\mu$ are integrally closed in $D$,
then $c = t$ works.
In particular,
if $D$ is a regular local ring, $c = t$ works.
\end{prop}

\proof
The variable $x$ induces a natural grading on $R$.

{\it Step one:}
We prove that there exists a positive integer $c$ such that
$\overline {x^tR+ \mu^nR} \subseteq x^tR+\mu^{\floor{n/c}}R$
when $n$ is a multiple of $t$.

As the ideal $x^tR+ \mu^nR$ is graded,
so is its integral closure.
Let $ax^l$ be a homogeneous element of this integral closure
with $a \in D$.
We may assume that $l < t$.
We write the equation of integral dependence:
$$
(ax^l)^k + b_1 (ax^l)^{k-1} + b_2 (ax^l)^{k-2} + \cdots + b_k = 0
$$
for some homogeneous $b_i \in \left(x^tR + \mu^nR\right)^i$.
We may assume that each summand in the equation above has $x$-degree
precisely $lk$.
Thus $b_i \in \left(x^tR + \mu^nR\right)^i \cap (x)^{li}$.
Write $b_i = x^{li} a_i$ for some
$a_i \in D \cap ((x^tR + \mu^nR)^i : x^{li})$.
As all these ideals are graded,
$$
D\cap ((x^tR + \mu^nR)^i : x^{li})
=
\sum _{j=0} ^i x^{tj} \mu^{n(i-j)}R
:\hskip -0.3em\raisebox{-.4ex}{$ \atop D$} x^{li}
=
\sum _{tj \le li} \mu^{n(i-j)}D
\subseteq \mu^{ni/t},
$$
where the last line follows from $tj \le li$ and $l < t$.
Thus after dividing the integral equation for $ax^l$ above by $x^{lk}$,
we see that $a$ is integral over $\mu^{n/t}$.
As $D$ is analytically unramified,
by Proposition \ref{prop:selfrad}\  
there exists an integer $k$ such that $\overline{\mu^n} \subseteq \mu^{n-k}$
for all $n$.
Thus there exists an integer $c \ge t$ ($c = t(k+1)$ will do)
such that
$\overline{\mu^n} \subseteq \mu^{\floor{n/c}}$
for all $n$.
This proves that
$ax^l \in x^t R + \mu^{\floor{n/c}}R$ whenever $n$ is divisible by $t$.

Note that $c = t$ works if the powers of $\mu$ are integrally closed.

{\it Step two:}
We prove that for the $c$ from Step~1,
$\overline {x^tR+ \mu^nR} \subseteq x^tR+\mu^{\floor{n/c}}R$
for all $n$.

Write $n = qc + r$ for some integers $q$ and $r$ with $0 \le r < c$.
Then
\begin{eqnarray*}
\overline {x^tR+ \mu^nR}
&\subseteq&
\overline {x^tR+ \mu^{qc}R} \\
&\subseteq&
x^tR+\mu^qR \mbox{\ \ \ by Step 1} \\
&=&
x^tR+\mu^{\floor{n/c}}R.
\end{eqnarray*}

{\it Step three:}
We prove that for the $c$ from Step~1,
$\overline {x^tR+ m^nR} \subseteq x^tR+m^{\floor{n/c}R}$
for all $n$.

If $c > n$,
$\floor{n/c} = 0$
so $m^{\floor{n/c}} = R$ and the inclusion holds trivially.
So we now assume that $t \le c \leq n$.
Then
$\overline {x^tR+ m^nR}
= \overline {x^tR+ x^nR + \mu^nR}
= \overline {x^tR+ \mu^nR}$.
By the previous step then
$\overline {x^tR+ m^nR}
\subseteq x^tR+\mu^{\floor{n/c}}R
\subseteq x^tR+m^{\floor{n/c}}R$,
so done.
\qed
\medskip

A similar proof works
if $x^t$ is replaced by a monomial in several variables.
However,
the following useful lemma enables an alternate proof to be given
in Corollary \ref{cor:prodmonomial}:

\begin{lemma}
\label{lemma:primdec}
Assume $I=J\cap K$, where $J$ and $K$ are integrally closed
ideals. If the
theorem holds for $J$ and $K$ then it holds for $I$.
\end{lemma}

\proof
By assumption $\overline {J+m^{n}}\subseteq J+m^{\floor {n/c}}$
and $\overline {K+m^{n}}\subseteq K+m^{\floor {n/c}}$.
Let $z\in \overline{I + m^n}$.
Then $z \in (J+m^{\floor {n/c}})\cap (K+m^{\floor {n/c}})$.
Write $z=j+m_{1}=k+m_{2}$, where $m_{1}, m_{2}\in m^{\floor {n/c}}$.
By the Artin-Rees Lemma
\begin{eqnarray*}
j-k\in m^{\floor {n/c}}\cap (J+ K)=
m^{\floor {n/c}-c_{1}}(m^{c_{1}}\cap (J+ K))
&\subseteq & m^{\floor {n/c}-c_{1}}(J+ K)\\
&\subseteq & m^{\floor {n/c}-c_{1}}J+m^{\floor {{n/c}-c_{1}}}K \\
&\subseteq & m^{\floor {n/c}-c_{1}}\cap J+m^{\floor {n/c}-c_{1}}\cap K.
\end{eqnarray*}
Then $j-k=x+y$, where $x\in m^{\floor {n/c}-c_{1}}\cap J$ and $y\in
m^{\floor {n/c}-c_{1}}\cap K$, so
$j-x\in J$, $k+y\in K$ and $j-x=k+y$.
Finally $j=k+y+x\in J\cap K +m^{\floor {n/c}-c_{1}}$. \qed

\begin{cor}
\label{cor:prdecpr}
Let $fR$ be a principal ideal in a local normal domain $(R,m)$.
If the theorem holds for each primary component of $fR$
then it holds for $fR$.
\end{cor}

\proof
Consider a primary decomposition of $fR$: $fR=Q_{1}\cap \cdots \cap Q_{r}$.
Since $R$ is a normal domain all the $Q_{i}$ have height 1.
Moreover, $Q_{i}=(fR)_{\sqrt {Q_{i}}}\cap R$
so each $Q_{i}$ is integrally closed as $(fR)_{\sqrt {Q_{i}}}$ is.
We now apply Lemma \ref{lemma:primdec}. \qed
\bigskip

In this corollary we used the fact that
the principal ideals are integrally closed in normal domains.
Now by applying Corollary \ref{cor:prdecpr} and Proposition \ref{prop:monomial}
we obtain the following:

\begin{cor}
\label{cor:prodmonomial}
Let $D$ be an analytically unramified local ring and
$\x$ variables over $D$. Let $R=D[[\x ]]$
and let $m$ be the maximal ideal.
Let $I=(x_1^{t_1}\cdots x_v^{t_v})R$.
Then $\overline {I+ m^n} \subseteq I+m^{\floor{n/c}}$
for some positive integer $c$ and for all $n$.
\qed
\end{cor}

Also,
if $I$ is generated by powers of some of the variables,
the main theorem also holds in $R$ for $I$.
The proof uses linear algebra arguments
and is given in Section~4.

\begin{prop}
\label{prop:propwei}
Let $(D,m_{D})$ be a complete local ring,
$\x $ variables over $D$, $R=D[[\x ]]$,
$f\in R$.
Assume that there exists an integer $h$ such that:\\
$i)$ the coefficients of $x^{i}_{v}$ in $f$ are in $(m_D, x_1,\ldots ,x_{v-1})$
for $0\leq i<h$;\\
$ii)$ the coefficient of $x^{h}_{v}$ in $f$ is not in $(m_D, x_1,\ldots ,x_{v-1})$.\\
Then the main theorem holds for $fR$.
\end{prop}

\proof
We can apply Weierstrass Preparation Theorem and write $f=uf^{*}$,
where $u$ is a unit in $D[[x_{1},\ldots ,x_{v-1}]]$
and $f^{*}$ is a monic polynomial in $D[[x_{1},\ldots ,x_{v-1}]][x_{v}]$.
Without loss of generality we may replace $f$ by $f^{*}$.
Let $A=D[[x_{1},\ldots ,x_{v-1}]]$. Let $A'$ be a finite (local)
extension of $A$ such that in $S=A'[[x_{v}]]$,
$f$ factors into linear factors, say:
$f=(x_{v}-\alpha _{1})^{r_{1}}\cdots (x_{v}-\alpha _{s})^{r_{s}}$.
By Proposition \ref{prop:monomial}, $(x_{v}-\alpha _{i})^{r_{i}}S+m^{n}_{S}\subseteq
(x_{v}-\alpha _{i})^{r_{i}}S+m^{\floor{n/{c'}}}_{S}$
for all $n>>0$.
By Corollary \ref{cor:prodmonomial},
$\overline {fS+m^{n}_{S}}\subseteq fS+m^{\floor{n/c}}_{S}$
for all $n>>0$.
By Lemma \ref{lemma:modulefinite} then the theorem holds in $R$. \qed
\medskip

\begin{cor}
Let $R$ be a complete regular local ring containing a field,
and $I=fR$ a principal ideal.
Then $\overline{I + m^n} \subseteq I + m^{\floor{n/c}}$
for some $c$ independent of $n$.
\end{cor}

\proof
By the Cohen Structure Theorem $R$ is a power series ring
$k[[x_{1},\ldots ,x_{q}]]$ over a field $k$.
If $k$ is finite,
let $S=\overline {k}[[x_{1},\ldots ,x_{q}]]$ where $\overline {k}$
is an algebraic closure of $k$.
If the theorem holds in S for $fS$, then
by faithful flatness it holds in $R$,
so without loss of generality we may assume that
$R$ contains an infinite field.
Write $f=f_{s}+\cdots$,
where the degree of $f_{s}$ is positive and lowest.
We can choose elements $u_{1},\ldots ,u_{q}$ in $k$ such that,
after the change of variables $y_{1}=x_{1},y_{2}=x_{2}-u_{1}x_{1},\ldots ,
y_{q}=x_{q}-u_{q}x_{1}$, $f_{s}$ is monic in $y_1$.
We then apply Proposition \ref{prop:propwei}. \qed

We need similar results for rings not containing fields:

\begin{prop}
\label{prop:dvr}
Let $V$ be a discrete valuation ring with maximal ideal $pV$.
Set $R=V[[\xx ]]$.
Then $\overline {p^{t}R+(\xx )^{n}R}\subseteq p^{t}R+(\xx )^{\floor {n/t}}R$.
\end{prop}

\proof
The variables $\xx $ induce a multigrading on $R$.
The ideal $p^{t}R+(\xx )^{n}R$ is multihomogeneous,
therefore its integral closure is homogeneous as well.
Let $a \in \overline {p^{t}R+(\xx )^{n}R}$.
We may assume that $a$ is multihomogeneous,
so without loss of generality
$a=p^{l}\ux ^{\nnu}$,
$l\leq t-1$,
$|\nnu|=\nu _{1}+\cdots +\nu _{d}\leq n-1$.
We write an equation of integral dependence of $a$
over $p^{t}R+(\xx )^{n}R$:
$$
(p^{l}\ux ^{\nnu})^{k}+b_{1}(p^{l}\ux ^{\nnu})^{k-1}+\ldots +b_{k}=0,
$$
with $b{_i}\in (p^{t}R+(\xx )^{n}R)^{i}$.
We may assume that each $b_{i}$ is homogeneous
and that the $x_{j}$-degree of the i$^{th}$ summand
is exactly $\nu _{j}(k-i)$.
Then $b_{i}=c_{i}\ux ^{i \underline {\nu}}$
with $c_i \in V$.
The integral equation now has the form:
$$p^{lk}\ux ^{\nnu k}+c_{1}p^{l(k-1)}\ux ^{\nnu k}+\cdots +c_{k}\ux ^{\nnu k}=0,$$
so for all $i$ we can write: $c_{i}p^{l(k-i)}=u_{i}p^{f(i)}$, where
$u_{i}$ is a unit in $V$.
Necessarily there exists at least one $i$ such that $f(i)\leq lk$.
Then
$$
a^{k}=p^{lk}\ux ^{\nnu k}
	=c_{i}\frac{p^{lk-f(i)}}{u_{i}}p^{l(k-i)}\ux ^{\nnu k}
	=\frac{p^{lk-f(i)}}{u_{i}} b_{i}a^{k-i}.
$$
So by possibly modifying the original integral equation we get
a homogeneous integral equation of the form $a^k=b_k$,
with $b_k\in (p^{t}R+(\xx )^{n}R)^k$.
As the equation is homogeneous,
there exists $i$, $1 \le i \le k$,
such that
$|\nnu |k\geq n(k-i)$ and $lk\geq ti$.
Assume for contradiction that $|\nnu |<\floor {n/t}\leq n/t$.
Then:
$$
\frac{nk}{t} > |\nnu| k \geq n(k-i)
\Longrightarrow
ti > (t-1)k\geq lk \ge ti,
$$
which is a contradiction.
\qed

\bigskip

Our goal is to prove the same in greater generality:
namely that for any ideal $fR$ in a complete normal local domain $(R,m)$,
any element $u$ of $\overline{fR + m^n}$ lies in $fR + m^{\floor{n/c}}$
for some $c$ independent of $u$ and $n$.
However, in this generality,
$c$ is not as easy to determine.

Since the proof of the general result is quite involved, we outline it here.

We first rephrase the statement:
we see that $Z = u$ satisfies an equation
$$
Z^d+Z_1 f Z^{d-1} + Z_2 f^2 Z^{d-2} + \cdots + Z_d f^d = 0,
$$
modulo $m^n$ for some $Z_1, \ldots, Z_d$.
If we could use the \LAA theorem,
we would be able to conclude that there exists a true solution of this equation
which differs from the original approximate solution by $m^{\floor{n/c}}$
for some $c$ independent of $n$
and hence that $u \in fR + m^{\floor{n/c}}$.
However,
there are two problems with this:
one is that the $d$ in the equation depends on $u$ and $n$,
and the other one is that we do not know of a proof of the general \LAA Theorem.
We overcome both difficulties:
in Theorem~\ref{theorem:bound}
we prove that there exists one $d$ which works for all $u$ and all $n$,
except that the resulting equation of integral dependence is not over
$fR + m^n$ but over $fR + m^{\floor{n/d}}$.
Also,
in Corollary~\ref{cor:pf} we prove a special case of the \LAA Theorem.
However,
for these intermediate results we need slight assumptions on the $f$
in case $R$ is a ring of mixed characteristic.
%
The proof of the general case in Theorem~\ref{theorem:general}
also uses the fact that all the primary components of
principal ideals are the symbolic powers of height one prime ideals,
and that such symbolic powers
are associated to powers of radical principal ideals
after by passing to a faithfully flat extension $S$ of $R$.

\medskip


\begin{theorem}
\label{theorem:bound}
Let $(R,m)$ be a complete normal local domain
and $fR$ a non-zero principal ideal.
In the case when $R$ does not contain a field,
we let
$p$ be a generator of the maximal ideal in a coefficient ring for $R$,
and we assume that $f$ satisfies one of the following properties:
(i) $f,p$ is a part of a system of parameters,
or (ii) $f=ap^c$ for some positive integer $c$
and some element $a$ of $R$
not contained in any minimal prime ideal over $pR$.

Then there exist integers $d$ and $l$ such that for each $n$,
every element in $\overline{fR+m^n}$
satisfies an integral equation of degree $d$ over $fR+m^{\floor {n/l}}$.
\end{theorem}

\proof
It is sufficient to prove that
if $J$ is $m$-primary,
then there exists an integer $d$ such that for each $n$,
every element in $\overline{fR+J^n}$ satisfies an integral equation of degree $d$
over $fR+J^{\floor {n/d}}$.
(Note, however, that $d$ depends on $J$!)

We use the Cohen Structure Theorem.
Let $f_1, \ldots, f_l$ be a system of parameters in $R$.
When $R$ contains a coefficient field $k$,
we may assume that $f_1 = f$,
and we define $A = k[[f_1, \ldots, f_l]]$.
When $R$ contains a coefficient ring $(V, (p))$ of dimension $1$,
we may assume that $p$ is $f_1$.
In case (i) we may also assume that $f = f_2$
and in case (ii) we may assume that $f_2$ is $a$ if $a$ is not a unit.
In case (ii) if $a$ is a unit,
as $fR = p^cR$,
without loss of generality $a = 1$.
We then define $A = V[[f_2, \ldots, f_l]]$.
In either case,
set $J = (f_1, \ldots, f_l) A$.
By the Cohen Structure Theorem,
$A$ is a regular local ring contained in $R$,
$R$ is module-finite over $A$,
and $JR$ is $m$-primary.
We will prove the theorem for this $JR$.
Furthermore,
we will prove that the integral equation of degree $d$
will have coefficients in $A$.

Let $K$ be the fraction field of $A$ and $L$ the fraction field of $R$.
By elementary field theory
there exist fields $L'$ and $F$
such that all the inclusions
$K \subseteq F \subseteq L'$ and $L \subseteq L'$ are finite,
such that $L'$ is Galois over $F$
and such that $F$ is purely inseparable over $K$.
To simplify notation,
as the coefficients of the integral equation will actually lie in $A$,
we may replace $R$ by the integral closure of $R$ in $L'$
and so we may assume that $L = L'$.
Let $d = [L : F]$ and $e = [F : K]$.
Let $S$ be the integral closure of $A$ in $F$.
Then $S$ is a complete normal local domain between $A$ and $R$
and the extension from $S$ to $R$ is Galois.

Let $u\in \overline {fR+(JR)^{n}}$.
Consider the (at most) $d$ conjugates of $u$ over $S$,
say $u = u_1, u_2, \ldots, u_d$.
Write an integral equation for $u$ over $fR+(JR)^n$:
$$
u^k + \alpha_1 u^{k-1} + \alpha_2 u^{k-2} + \cdots + \alpha_k = 0
$$
with $\alpha _{i}\in (fR+(JR)^{n})^{i}$.
By applying field automorphisms to this equation
and by using that $(f)$ and $J$ are ideals of $A$ (and thus of $S$),
we obtain that each $u_i$ is integral over $fR + J^nR$.
Let $s_{h}$ be the sum of the products of the $u_{i}$,
taken $h$ at a time ($h$th symmetric function in the $u_i$).
Then
$$
u^d - s_1 u^{d-1} + \cdots + (-1)^d s_d = 0,
$$
and
$s_{h} \in \overline {(fR+(JR)^{n})^{h}}\cap S$.
We raise all this to the $e$th power.
As $e$ is either $1$ or a power of the characteristic $p$ of the given fields,
we obtain
\begin{eqnarray*}
u^{de} &-& s_1^e u^{e(d-1)} + \cdots + (-1)^{de} s_d^e = 0, \mbox{\ \ and} \\
s_h^e &\in& \overline {(fR+(JR)^n)^h}^e\cap A \\
&\subseteq&
\overline {(f^{he}R+(JR)^{nhe})} \cap A \\
&\subseteq&
\overline {(f^{he}A+(JA)^{nhe})}, \\
\end{eqnarray*}
as $A \subseteq S$ is a module-finite extension.
By Propositions~\ref{prop:monomial} and~\ref{prop:dvr}
and by Corollary~\ref{cor:prdecpr} then there exists an integer $l$
such that
$$
s_h^e \in f^{he}A+(JA)^{\floor{n/l}}
\subseteq
(fA+(JA)^{\floor{n/lhe}})^{he}
\subseteq
(fA+(JA)^{\floor{n/lde}})^{he}.
$$
Thus $u$ satisfies an equation of integral dependence of degree $de$
over $fR + (JR)^{\floor{n/lde}}$,
all of whose coefficients are in $A$.
\qed
\medskip


Now we are in the set-up of the \LAA theorem (from two pages back).
We have an equation
$$
Z^d+Z_1 f Z^{d-1} + Z_2 f^2 Z^{d-2} + \cdots + Z_d f^d = 0, \,\,\,
$$
where $Z$ and the $Z_{i}$ are indeterminates,
which is independent of $n$ and $u$.
We are given a solution
$(Z = u, Z_1 = z_1, \ldots, Z_d = z_d)$
of such an equation modulo a high power $m^n$ of the
maximal ideal of $R$.
The following theorem proves that under some assumptions on $f$,
we can find a true solution
$Z = w,Z_1 = w_1, \ldots, Z_d = w_d$
of this equation such that $w - u \in m^{\floor{n/c}}$
for some $c$ independent of $u$ and $n$.
This proves that $u \in fR + J^{\floor{n/c}}R$,
and it also proves a special case of the \LAA Theorem:

\begin{prop}
\label{prop:power}
Let $(R,m)$ be a Noetherian local integrally closed integral domain,
and $f \in R$  satisfying the following:
\begin{enumerate}
\item
There exists a positive integer $c$
such that for all $n \ge 1$,
$\overline{(f) + m^n} \subseteq (f) + m^{\floor{n/c}}$.
\item
For every $k = 1, \ldots, N$
there exist positive integers $d$ and $l$ such that for all $n$,
every element of 
$\overline{(f^k) + m^n}$ satisfies an equation of integral dependence
of degree $d$ over $(f) + m^{\floor{n/l}}$.
\end{enumerate}
Then for every $k = 1, \ldots, N$,
there exists a positive integer $c$ such that
$\overline{(f^k) + m^n} \subseteq (f^k) + m^{\floor{n/c}}$.
\end{prop}

\proof
We prove this by induction on $k$.
The case $k=1$ is assumed.
So assume $k>1$.
By induction,
$\overline {(f^{k})+m^{n}}\subseteq (f^{k-1})+m^{\floor {n/{c'}}}$
for some constant $c'$ independent of $n$.
We pick an element $u$ in $\overline {(f^{k})+m^{n}}$.
Write $u = rf^{k-1} + s$ for some $r \in R$ and $s \in m^{\floor{n/c'}}$.
It suffices to prove that $rf^{k-1}$ lies in $(f^k) + m^{\floor{n/c}}$
for some $c$ independent of $n$ and $u$.
Note that $rf^{k-1}$ is integral over $(f^k) + m^{\floor{n/c'}}$.
Hence it suffices to prove that
$(f^{k-1}) \cap \overline{(f^k)+m^{\floor {n/{c'}}}}$
is contained in $(f^k) + m^{\floor{n/c}}$ for some $c$ independent of $n$,
or even that
$(f^{k-1}) \cap \overline{(f^k)+m^n}$
is contained in $(f^k) + m^{\floor{n/c}}$ for some $c$ independent of $n$.
Thus without loss of generality we may assume that $u = rf^{k-1}$.
Our goal is to prove that $r \in \overline{(f) + m^{\floor{n/c''}}}$
for some integer $c''$ independent of $n$ and $r$,
for then we know that $r \in (f) + m^{\floor{n/c'''}}$ for some
$c'''$ independent of $n$ and $r$,
which proves that $u$ lies in the desired ideal.

We first prove that a power of $r$ lies in a good ideal,
and for that we need the following detour:
\smallskip

{\it Claim:}
$r^{d}\in (f)+m^{\floor {n/l}-e}$ for some constant $e$ independent of $n$.

\smallskip

{\it {Proof of the claim:}}
By assumption there exists an integer $d$ independent of $n, u$ and $r$
such that $rf^{k-1}$ satisfies an integral equation of degree $d$
over $(f^{k})+m^{\floor {n/l}}$, say:
$(rf^{k-1})^d + \alpha_1 (rf^{k-1})^{d-1}+\cdots +\alpha _{d}=0$,
where $\alpha _{i}\in ((f^{k})+m^{\floor {n/l}})^{i}$.

We will recursively define
$\beta _{d-i+1} \in ((f^k)+m^{\floor {n/l}})^{d-i}$
for each $i\in \{ 0,\ldots ,d-1\}$
such that
$$
r^{d}(f^{k-1})^{d-i} +\alpha _{1} r^{d-1} (f^{k-1})^{d-i-1} + \cdots
+\alpha _{d-i} r^{i} +\beta _{d-i+1} =0.
\ \ \ \ (\#)
$$
If $i=0$, set $\beta _{d+1}=0$.
Now assume we have defined $\beta _{d-i+1}$ for some $i < d-1$.
By the Artin-Rees Lemma there exists a positive integer $e$ such that
$m^n \cap (f^{k-1}) \subseteq f^{k-1} m^{n-e}$
for all $n \ge e$.
In the following we may and do assume that $n/l \ge e$.
With this we construct the next $\beta$ using the equation displayed above
and the following:
\begin{eqnarray*}
\alpha _{d-i}r^{i}+\beta _{d-i+1}
&\in& (f^{k-1}) \cap ((f^{k})+m^{\floor{n/l}})^{d-i} \\
&=& (f^{k-1}) \cap
	\left( f^{k} ((f^{k})+m^{\floor {n/l}})^{d-i-1}
		+ m^{\floor {n/l}(d-i)} \right) \\
&=& f^{k}((f^{k})+m^{\floor {n/l}})^{d-i-1}
	+ (f^{k-1}) \cap m^{\floor {n/l}(d-i)} \\
&\subseteq & f^{k-1}((f^{k})+m^{\floor {n/l}})^{d-i-1}
	+ (f^{k-1}) m^{\floor {n/l}(d-i) - e} \\
&\subseteq & f^{k-1}((f^{k})+m^{\floor {n/l}})^{d-i-1}
\end{eqnarray*}
as $n/l \ge e$.
Thus we may write $\alpha _{d-i}r^{i}+\beta _{d-i+1} = f^{k-1}\beta _{d-i}$
for some $\beta _{d-i}\in ((f^{k})+m^{\floor {n/l}})^{(d-i-1)}$.
To finish the induction step
we only have to divide the displayed equation (\#)
by the nonzerodivisor $f^{k-1}$.

In the final step $i = d-1$
we thus obtain $r^{d}f^{k-1}+\alpha _{1}r^{d-1}+\beta _{2}=0$.
Therefore
$$r^{d}f^{k-1}= -\alpha _{1}r^{d-1}-\beta _{2}
\in (f^{k-1})\cap ((f^{k})+m^{\floor {n/l}})
= (f^{k})+(f^{k-1})\cap m^{\floor {n/l}}
\subseteq (f^{k})+f^{k-1}m^{\floor {n/l}-e}.$$
It follows that $r^{d}\in (f)+m^{\floor {n/l}-e}$.
This completes the proof of the claim.

\bigskip
Now we are ready to prove that $r$ is integral over
$(f) + m^{\floor{n/dlk(e+1)}}$.
Recall that $rf^{k-1}\in \overline {(f^{k})+m^{n}}$.
It suffices to prove that
for any valuation $v$ on the field of fractions of $R$,
$v(r) \ge \min\{v(f), \floor{n/dlk(e+1)}v(m)\}$.

Since $rf^{k-1}\in \overline {(f^{k})+m^n}$,
$v(r) + (k-1) v(f) = v(rf^{k-1}) \ge \min \{kv(f),n v(m)\}$,
therefore $v(r) \ge \min \{v(f),n v(m) - (k-1)v(f)\}$.
If $v(r) \ge v(f)$, there is nothing to show,
so we may assume that
$$
\floor{n/(e+1)} v(m) - (k-1) v(f) \le n v(m) - (k-1)v(f) \le v(r) < v(f).
$$
This implies that $\floor{n/(e+1)} v(m) < kv(f)$.
Now we use our detour:
as $r^{d}$ lies in $(f)+m^{\floor {n/l}-e}$
$\subseteq (f)+m^{\floor {n/l(e+1)}}$,
then
$$
dv(r) \ge \min\{v(f),\floor{n/l(e+1)}v(m)\}
\ge \min\{v(f),\floor{n/lk(e+1)} v(m)\}.
$$
If
$dv(r) \ge \floor{n/lk(e+1)} v(m)$,
we are done,
so we may assume instead that
$$
\floor{n/lk(e+1)} v(m) > dv(r) \ge v(f).
$$
Thus
$$
\floor{n/lk(e+1)} v(m) > dv(r) \ge v(f) >
\frac{1}{k} \floor{n/(e+1)} v(m),
$$
which is a contradiction.
This finishes the theorem.
\qed
\bigskip

\begin{cor}
\label{cor:pf}
Let $(R,m)$ be a complete local normal domain
and let $(f)$ be a principal radical ideal.
In case $R$ does not contain a field,
let $(V,(p))$ be a general coefficient ring of $R$
and we also assume that either $fR = pR$
or that $f, p$ is part of a system of parameters in $R$.
Then for all $k$,
$\overline {(f^{k})+m^{n}}\subseteq (f^{k})+m^{\floor{n/c}}$
for some constant $c$ independent of $n$.
\end{cor}

\proof
The case $k=1$ holds by Proposition~\ref{prop:selfrad}.
Thus condition 1. of the previous theorem is satisfied.
Condition 2. of the previous theorem is satisfied
by Theorem~\ref{theorem:bound},
so that the corollary follows by the previous theorem,
Proposition~\ref{prop:power}.
\qed

Before we prove the main theorem, we need two more lemmas:

\begin{lemma}
\label{lemma:colon}
Let $(R,m)$ be a Noetherian local domain.
Suppose that for some ideal $I$ there exists an integer $c$
such that $\overline{I + m^n} \subseteq \oI + m^{\floor{n/c}}$
for all $n$.
Then for any nonzero element $y$ in $R$
there exists an integer $d$ such that
$\overline{I:y + m^n} \subseteq \oI:y + m^{\floor{n/d}}$.
\end{lemma}

\proof
By the Artin-Rees lemma there exists an integer $k$
such that in $R$ modulo $\oI$,
$(y)R + \oI/\oI \cap m^nR + \oI/\oI \subseteq y m^{n-k}R + \oI/\oI$.
Thus $(y) \cap (\oI + m^n) \subseteq y m^{n-k} + \oI \cap (y)$.

Let $u$ be in
$\overline{I:y + m^n}$.
Then $uy$ lies in
$\overline{I + m^n} \subseteq \oI + m^{\floor{n/c}}$
and also in $(y)$.
By the above then
$uy$ lies in
$\oI \cap (y) + y m^{\floor{n/c}-k}$,
so that
$u$ lies in
$\oI : (y) + m^{\floor{n/c}-k}$.
This finishes the lemma.
\qed

\begin{lemma}
\label{lemma:inteqn}
Let $R$ be an integral domain,
$x$ and $y$ non-zero elements of $R$
and $d, l$ positive integers
such that for every positive integer $n$,
every element of $\overline{(xy) + m^n}$
satisfies an integral equation of degree $d$ over
$(xy) + m^{\floor{n/l}}$.
Then there exists a positive integer $k$
such that for every positive integer $n$,
every element of $\overline{(x) + m^n}$
satisfies an integral equation of degree $d$ over
$(x) + m^{\floor{n/k}}$.
\end{lemma}

\proof
Let $r \in \overline{(x) + m^n}$.
Then $ry \in \overline{(xy) + m^n}$.
Thus there exist elements $r_i \in ((xy) + m^{\floor{n/l}})^i$
such that
$$
(ry)^d + r_1 (ry)^{d-1} + \cdots + r_{d-1} ry + r_d = 0.
$$
Write $r_i = s_i (xy)^i + t_i$
for some $s_i \in R$ and some $t_i \in m^{\floor{n/l}}$.
Then
$$
(ry)^d + s_1 (xy) (ry)^{d-1} + \cdots + s_{d-1} (xy)^{d-1} ry + s_d (xy)^d
+ t_1 (ry)^{d-1} + \cdots + t_{d-1} ry + t_d = 0.
$$
Thus
$t_1 (ry)^{d-1} + \cdots + t_{d-1} ry + t_d \in (y^d) \cap m^{\floor{n/l}}$.
By the Artin-Rees Lemma then there exists an integer $k$ such that
$t_1 (ry)^{d-1} + \cdots + t_{d-1} ry + t_d \in y^d m^{\floor{n/k}}$.
But then dividing the integral equation above by $y^d$
shows that $r$ satisfies an integral equation of degree $d$
over $(x) + m^{\floor{n/kd}}$.
\qed

With this we can prove the general result for principal ideals
in complete normal local domains:

\begin{theorem}
\label{theorem:general}
Let $(R,m)$ be a complete normal local domain.
Let $f$ be an element in $R$.
Then there exists a positive integer $c$ such that
$$
\overline{(f)+m^n} \subseteq
(f) +m^{{\floor{{n/c}}}}
{\mbox{\rm \ \ for all \ }} n.
$$
\end{theorem}

\proof
If $f = 0$,
the theorem is known by Proposition~\ref{prop:selfrad}.
So we may assume that $f$ is not zero.

As $R$ is normal,
all the associated prime ideals of the ideal $(f)$ are minimal over $(f)$.
By Corollary \ref{cor:prdecpr} it suffices to prove the theorem
for the primary components of $(f)$ in place of $(f)$.
Let $P$ be an associated prime ideal of $(f)$.
As $R$ is normal,
the localization $R_P$ is a one-dimensional regular local ring,
so $fR_P = P^kR_P$ for some integer $k$.
Thus the $P$-primary component of $fR$ equals
the $k$th symbolic power $P^{(k)}$ of $P$
and it suffices to prove the theorem for all $P^{(k)}$ in place of $(f)$.

Let $P = (a_1, \ldots, a_l)$.
Let $X_1, \ldots, X_l$ be indeterminates over $R$
and let $S$ be the faithfully flat extension
$R[X_1, \ldots, X_l]_{mR[X_1, \ldots, X_l]}$ of $R$.
Note that as all the associated primes of $xS$ have height one
and as $S$ localized at height one prime ideals is a principal ideal domain,
the ideal generated by $x = a_1X_1 + \cdots + a_l X_l$ is radical.

Suppose that this $x$ satisfies the conditions of Corollary \ref{cor:pf}.
Namely,
either $R$ contains a field,
or instead if $(V,(p))$ is a coefficient ring of $R$,
then either $x = p$ or $x, p$ is a part of a system of parameters.
Then by Corollary \ref{cor:pf} there exists an integer $c$ such that
$\overline{x^kS + m^nS} \subseteq x^kS + m^{\floor{n/c}}S$
for all $n$.
Note also that $PS$ is associated to $xS$
and that the $PS$-primary component of $x^kS$ is $P^{(k)}S$
(as $S_{PS}$ is a principal ideal domain).
Thus there exists an element $y$ in $S$ such that $x^kS : y = P^{(k)}S$.
As $R$ is normal,
then so is $S$,
so that $x^kS = \overline{x^k S}$.
An application of
Lemma~\ref{lemma:colon}
shows that
there exists an integer $c'$ such that
$\overline{P^{(k)}S + m^nS} \subseteq P^{(k)}S + m^{\floor{n/c'}}S$
for all $n$.
Finally,
\begin{eqnarray*}
\overline{P^{(k)} + m^n} &\subseteq& \overline{P^{(k)}S + m^nS} \cap R \\
&\subseteq& \left(P^{(k)}S + m^{\floor{n/c'}}S \right) \cap R \\
&=& P^{(k)} + m^{\floor{n/c'}}
\end{eqnarray*}
as $S$ is faithfully flat over $R$.

This finishes the theorem for rings containing fields.

Now assume that $R$ contains a coefficient field $(V, (p))$.
The above proves the theorem for all $f$
which are not contained in any minimal prime ideal over $pV$.
Thus by Lemma~\ref{lemma:colon},
for all height one prime ideals $P$ of $R$ not containing $p$
and all positive integers $k$
there exists an integer $c$ such that
$\overline{P^{(k)} + m^n} \subseteq P^{(k)} + m^{\floor{n/c}}$.

Let $P_1, \ldots, P_N$ be all the prime ideals in $R$
minimal over $pR$.
Let $W = R \setminus (P_1 \cup \cdots \cup P_N)$.
As $R$ is normal,
$W^{-1}R$ is a one-dimensional semi-local regular ring,
thus a principal ideal domain.
Let $x_i \in R$ such that $x_i W^{-1}R = P_i W^{-1} R$.
Therefore we may write
$p = u' x_1^{n_1} \cdots x_N^{n_N}$
for some unit $u' \in W^{-1}R$.
But then there exist $u, v \in W$ such that in $R$,
$up = v x_1^{n_1} \cdots x_N^{n_N}$.
Note that either $u$ is a unit in $R$
or else $p, u$ is a part of a system of parameters.
Thus by Theorem~\ref{theorem:bound},
for each positive integer $k$ there exist integers $d$ and $l$
such that every element of $\overline{(up)^k + m^n}$
satisfies an equation of integral dependence of degree $d$ over
$(up)^k + m^{\floor{n/l}}$.
Thus by Lemma~\ref{lemma:inteqn},
for each positive integer $k$ there exist integers $d$ and $l$
such that for all $i = 1, \ldots, N$,
every element of $\overline{(x_i)^k + m^n}$
satisfies an equation of integral dependence of degree $d$ over
$(x_i)^k + m^{\floor{n/l}}$.
This means that condition 2. of Proposition~\ref{prop:power}
is satisfied for each $x_i$.
But $x_iR = P_i \cap Q_i$,
where $Q_i$ is either the unit ideal
or a height one ideal modulo which $p$ is a non-zerodivisor.
As $P_i$ is a radical ideal (even prime),
by Proposition~\ref{prop:selfrad},
there exists a positive integer $c$ such that for all $n \ge 1$,
$\overline{P_i + m^n} \subseteq P_i + m^{\floor{n/c}}$.
By what we have proved,
there exists a positive integer $c'$ such that for all $n \ge 1$,
$\overline{Q_i + m^n} \subseteq Q_i + m^{\floor{n/c'}}$.
Thus by Lemma~\ref{lemma:primdec},
the theorem holds for $x_i$.
In particular,
condition 1. of Proposition~\ref{prop:power} is satisfied for $x_i$.
Thus by Proposition~\ref{prop:power},
theorem holds for all $x_i^k$,
as $k$ varies over all positive integers.
Then by Lemma~\ref{lemma:colon},
there exists an integer $c$ such that for all $i = 1, \ldots, N$,
$$
\overline{P_i^{(k)} + m^n} \subseteq P_i^{(k)} + m^{\floor{n/c}}.
$$
Thus by Lemma~\ref{lemma:primdec},
the theorem holds for $f = p$.

Hence condition 1. of Proposition~\ref{prop:power} is satisfied for $p$,
and condition 2. is satisfied by Theorem~\ref{theorem:bound}.
Thus by Proposition~\ref{prop:power},
the theorem holds for each $f = p^k$.

It remains to examine the case when $f$ and $p$
do not form a system of parameters.
In this case there exist an integer $e$,
an element $u \in W$
and $h \in R$,
such that $fh = up^k$.
We know the theorem for $uR$ and $p^kR$.
Since $uR$ and $p^kR$ are part of a system of parameters,
by Lemma \ref{lemma:primdec}
we also know the theorem for $(u) \cap (p^e) = (up^e) = (fh)$.
This means that there exists an integer $c$ such that
$\overline {fhR+m^{n}} \subseteq fhR + m^{\floor{n/c}}$.

Now pick $u\in \overline {fR+m^{n}}$.
Then $hu\in \overline {fhR+m^{n}} \subseteq fhR + m^{\floor{n/c}}$,
so
$hu \in fhR + m^{\floor{n/c}} \cap hR$.
By the Artin-Rees Lemma there exists an integer $k$ independent of $u$ and $n$
such that
$m^{\floor{n/c}} \cap hR \subseteq h m^{\floor{n/c}-k}$.
Thus
$hu \in fhR + h m^{\floor{n/c} - k}$,
so
$u \in fR + m^{\floor{n/c} - k}$.
This finishes the proof of this theorem.
\qed
\medskip

By using the reductions from Section~2,
the last theorem now
completely proves the main theorem stated in the introduction.
We restate it here for completeness.

\begin{theorem}
{\bf {\em {(The Main Theorem)}}}
Let $(R,m)$ be an excellent local ring.
Let $I$ be an ideal of $R$.
Then there exists a positive integer $c$ such that
$$
\overline{I+m^n} \subseteq
\oI +m^{{\floor{{n/c}}}}
{\mbox{\rm \ \ for all \ }} n.
$$
\end{theorem}

\section{Explicit bounds}

In the following examples we calculate explicitly the constant
$c$ from the main theorem.

\bigskip
The first ring we consider is an isolated hypersurface singularity,
for which the \LAA Theorem holds.

\begin{example}
\label{example:const}
Let $k$ be a field, $a$ a positive integer, $R=k[[X,Y,Z]]/{(X^{a}+Y^{a}+Z^{a})}=
k[[x,y,z]]$, $I=x^{t}R$.
To avoid the trivial case,
we assume that $a \ge 2$.
If char $k$ is positive, say $p=$ char $k$,
and $a=bp^{e}$ for some $b$ and $e$,
then
$X^{a}+Y^{a}+Z^{a}=(X^{b}+Y^{b}+Z^{b})^{p^{e}}$.
Since by Lemma \ref{lemma:reduced}
it suffices to find the bound in $R_{red}=k[[X,Y,Z]]/{(X^{b}+Y^{b}+Z^{b})}$,
thus without loss of generality we may assume that $e = 0$
and that char $k$ does not divide $a$.

Let $A=k[[x^{t},y]]\subseteq R$.
Since we are assuming that char $k$ does not divide $a$,
the fraction field of $R$,
$ff\, (R)$ for short,
is a separable field extension of $ff\,(A)$.
If necessary,
we adjoin the $a$-th roots of $1$ to $k$,
getting a field $K$ such that the extension $K/k$ is Galois.
The extension $ff\,(K[[x,y,z]])/ff\,(A)$ is also Galois.
As $[K:k]=\Phi (a)$, with $\Phi (a)$ a factor of the Euler function of $a$,
we have:
$[ff\,(K[[x,y,z]]):ff\,(R)]=\Phi (a)$.
Since $[ff\,(R):ff\,(k[[x,y]])]=a$
and $[ff\,(k[[x,y]]):ff\,(A) ]=t$,
we have that $[ff\,(K[[x,y,z]]):ff\,(A)]=ta\Phi (a)$.

Set $S=K[[x,y,z]]$ and $d=ta\Phi (a)$.
We pick $u\in \overline {x^{t}S+(x,y,z)^nS}$.
We have shown that the equation
$T^d+T_1 x^{t} T^{d-1} + T_2 x^{2t} T^{d-2} + \cdots + T_d x^{td} = 0$
has a solution $(u,z_{1},\ldots ,z_{d})$ mod $y^{\floor {n/d}}A$
(see comment before Proposition \ref{prop:power}).
Set
$\mu = \rank _{k} \frac {K[[X,Y,Z]]}{(X^{a-1},Y^{a-1},Z^{a-1})}=(a-1)^{3}$
(rank of the ring modulo the Jacobian ideal of $(X^{a}+Y^{a}+Z^{a})$).
By
\cite[\S 3]{Ljlaa},
the equation has a true solution, say $(w,w_{1},\ldots ,w_{d})$,
such that
$$
ord(u-w)\geq \frac{\floor {n/d}-a-1}{a\mu }\geq
\floor {\frac {n}{da^{2}\mu }} = \floor {\frac {n}{da^{2}(a-1)^{3} }}
$$
and so
$u-w \in (x,y,z)^{\floor {n/{da^{2}(a-1)^{3} }}}S$.
We conclude that $u\in x^{t}S+(x,y,z)^{\floor {n/{da^{2}(a-1)^{3} }}}S$.
Since $S=K[[x,y,z]]$ is faithfully flat over $R=k[[x,y,z]]$,
we have that $\overline {x^{t}R+m^{n}}\subseteq
x^{t}R+m^{\floor {n/{da^{2}(a-1)^{3} }}}$.

If $t=1$ then we do not need \LAA Theorem.
We can find an effective bound just by using
Proposition \ref{prop:selfrad}.
It is easy to see that the ideal $(x)+m^{na}$ is integrally closed. Then
$\overline {(x)+m^{n}}\subseteq (x)+m^{\floor {n/a}a}
\subseteq (x)+m^{n-(a-1)}$.
\end{example}

\bigskip
The next example is a generalization of Proposition \ref{prop:monomial}:

\begin{example}
Let $k$ be a field,
$X_1, \ldots, X_d$ variables over $k$ and $R=k[[X_{1},\ldots ,X_{d}]]$.
Set $I = (X_1 ^{a_1}, \ldots, X_s^{a_s})$
with each $a_i$ a positive integer.
Then for all $n$,
$$
\overline{I+m^n} \subseteq \oI + m^{\floor{{n \over a_1 \cdots a_s}}}.
$$

\proof
(Some of the reasoning below applies to arbitrary monomial ideals.)
We write a monomial $u=X_1^{b_1} \cdots X_d^{b_d}$ as a vector
$\vec \beta =(b_1, \ldots, b_d)$.
It is well known and easy to verify that a monomial $u$ is in the integral
closure of a monomial ideal $J$ if and only if
the vector $\vec \beta$ corresponding to $u$
is in the ``infinite'' convex hull of the vectors
$\vec \alpha_1, \ldots, \vec \alpha_r$
corresponding to the generators of $J$,
i.e.\ if and only if
$$
\vec \beta \ge \sum_{i=1}^r t_i \vec \alpha_i
\mbox{\ \ componentwise for some $t_i \ge 0$ satisfying \ \ }
\sum_{i=1}^r t_i = 1
\ \ \ \ \ \ \ \ (\star \star)
$$
In our case $J = I + m^n$
and
$\overline{I + m^n} =
\overline{(X_1 ^{a_1}, \ldots, X_s^{a_s},X_{s+1}^n, \ldots, X_d^n)}$
for $n \ge \max\{a_1, \ldots, a_s\}$.
Thus in our case we may assume that $r\le d$ and that
all the $\vec \alpha_j$ are of the form $c_j \vec e_j$,
where the vector $\vec e_j$ has $1$ in the $j$th row and $0$ everywhere else
and $c_j$ is either $n$ or $a_i$.
Without loss of generality we assume that $c_i = a_i$ for $i = 1, \ldots, p$
and $c_i = n$ for $i > p$.

We want to prove that if $u = X_1^{b_1} \cdots X_d^{b_d}$
is integral over $I+ m^n$,
then either $u \in \oI$ or else $u$
lies in $m^{\floor{\frac{n}{a_1 \cdots a_s}}}$.
In vector notation this
says that given $\vec \beta = (b_1, \ldots, b_d)$
satisfying ($\star \star$),
either
$\sum_{i=1}^s {b_i \over a_i} \ge 1$ or else
$|\vec \beta| = \sum_{i} b_i \ge n / a_1 \cdots a_s$.
We assume that
$\sum_{i=1}^s {b_i \over a_i} < 1$
and hence we have to prove that
$\sum_{i} b_i \ge n / a_1 \cdots a_s$.

We may rewrite ($\star \star$) in matrix notation as
$A \vec y = \vec \gamma$,
where
$$
A = \left[ \begin{array}{cccc}
\vec \alpha_1 & \cdots & \vec \alpha_r & \vec e_j \\
1 & \cdots & 1 & 0 \\
\end{array}
\right], \mbox{\hskip 0.5cm}
\vec y = \left[ \begin{array}{c}
t_1 \\
\vdots \\
t_r \\
q_j \\
\end{array}
\right], \mbox{\hskip 0.5cm}
\vec \gamma = \left[ \begin{array}{c}
b_1 -q_1 \\
\vdots \\
b_j \\
\vdots \\
b_d -q_d \\
1 \\
\end{array}
\right]
$$
and the $q_i$ are some nonnegative numbers.
By removing any of the $\vec \alpha_i$ if necessary
we may assume that ($\star \star$)
has all the $t_i$ strictly positive
and moreover that it cannot be written with any $t_i$ being zero.
As we want to prove that
$|\beta| \ge n / a_1 \cdots a_s$,
it suffices to prove the same
after omitting any of the first $d$ rows in $A$,
$\vec \beta$ and $\vec \gamma$.
Thus we may assume that $r = d$.

The first step is to make $d-1$ of the $q_i$ zero.
Note that $A$ is an invertible matrix.
Thus by applying Cramer's rule,
we see that the entries in $\vec y$
depend linearly on the entries of $\vec \gamma$.
We now decrease all the positive $q_i$'s, $i\not = j$:
while decreasing the $q_i$'s,
we want to keep all the $t_i$ and all the $q_i$ nonnegative.
If the given $q_j$ ever becomes zero in this way,
we exchange the roles of this $q_j$ with some nonzero $q_i$
(and also modify the last column of $A$),
and repeat again.
Note that by the choice of the $\alpha_i$,
none of the $t_i$ can become zero.
In this way, we obtain $d-1$ of the $q_i$ to be zero
and $q_j \ge 0$.

The cases $j \le p$ and $j > p$ differ only slightly
so we only finish the proof in the case $j = d$.
Cramer's rule gives
$$
q_d = \sum_{i=1}^d {b_i \over c_i} n - n,
t_i = {b_i \over a_i} \mbox{\ for $i \le p$},
t_i = {b_i \over n} \mbox{\ for $i = p+1, \ldots, d-1$}
$$
(i.e., for $i < d$, $t_i = {b_i \over c_i}$).
Now $q_d \ge 0$ says that
$\sum_{i=p+1}^d b_i \ge n(1 - \sum_{i=1}^p {b_i \over a_i})$.
We assumed $1 - \sum_{i=1}^p {b_i \over a_i} > 0$.
Thus
$1 - \sum_{i=1}^p {b_i \over a_i}$
equals a positive integer divided by $a_1 \cdots a_p$,
hence
$|\vec \beta| \ge \sum_{i=p+1}^d b_i \ge {n \over a_1 \cdots a_p}$,
and finally
$|\vec \beta| \ge {n \over a_1 \cdots a_s}$.
\qed
\end{example}

\vskip 2cm
\normalsize
\renewcommand{\baselinestretch}{0.7}

\vskip 1cm

\noindent
\sc
Department of Mathematics, University of Michigan, Ann Arbor, MI 48109-1003

\vskip 5mm
\noindent
Department of Mathematical Sciences, New Mexico State University, Las Cruces,
NM 88003-8001
\end{document}